\documentclass[10pt]{article}
\usepackage{amssymb}
\usepackage{enumerate}
\usepackage{graphicx}
\usepackage{amsthm}
\usepackage{color}
\usepackage{epstopdf}

\def\R{{\rm I\! R}}

\def\sss{strictly star-shaped }
\def\hkx{{\partial H_{2k} \over \partial x}}
\def\hky{{\partial H_{2k} \over \partial y}}

\newtheorem{theorem}{Theorem}

\newtheorem{lemma}{Lemma}

\title{Existence and uniqueness of limit cycles in a class of second order ODE's with inseparable mixed terms }

\author{M. Sabatini  
\footnote{Dip. di Matematica, Univ. di Trento, I-38050 Povo, (TN) - Italy.
Email: marco.sabatini@unitn.it,
Phone: ++39(0461)881670, Fax: ++39(0461)881624 - A previous version of this paper can be found  at www.arxiv.org as arXiv:1003.0803v1 [math.DS].}
}
\date{April $15t^h$, 2010}
\begin{document}
\maketitle
\begin{abstract}
We prove a uniqueness result for limit cycles of the second order ODE $\ddot x + \dot x \phi(x,\dot x) + g(x) = 0$. Under mild additional conditions, we show that such a limit cycle attracts every non-constant solution. As a special case, we prove limit cycle's uniqueness for an ODE studied in \cite{ETA} as a model of pedestrians' walk. 
This paper is an extension to equations with a non-linear $g(x)$ of the results presented in \cite{S}.

\end{abstract}

\section{Introduction}
 
The simplest non-linear continuous dynamical systems originate from the study of planar differential systems,
\begin{equation}\label{sysPQ}
\dot x = P(x,y), \qquad \dot y = Q(x,y), \qquad P, Q \in C^1(\R^2,\R^2).
\end{equation} 
Special cases of such systems are Lotka-Volterra ones, and systems equivalent to Li\'enard equations, 
\begin{equation}\label{syslie} 
\dot x = y - F(x)   , \qquad \dot y = - g(x),
\end{equation} 
or to Rayleigh equations
\begin{equation}\label{sysray} 
\dot x = y , \qquad \dot y = - g(x) - f(y).
\end{equation} 
All of them arise as mathematical models of biological, physical, engineering systems \cite{LC}. The study of the dynamics of (\ref{sysPQ}) strongly depends on the existence and stability properties of special solutions such as equilibrium points and non-constant isolated periodic solutions. In particular, if an attracting non-constant periodic solution exists, then it dominates the dynamics of (\ref{sysPQ})  in an open, connected subset of the plane, its region of attraction. Studying the number and location of isolated periodic solutions, usually called limit cycles, is by no means a trivial question, as shown by the resistance of Hilbert XVI problem (see \cite{Sm}, problem 13). In some cases such a region of attraction can even extend to cover the whole plane, with the unique exception of an equilibrium point. In such a case the limit cycle is unique and dominates the system's dynamics, as in \cite{DO}.
Uniqueness theorems for limit cycles have been  extensively studied (see \cite{CRV}, \cite{XZ1}, \cite{XZ2},  for recent results and extensive bibliographies). Limit cycle's uniqueness is a relevant feature even in discrete time systems, which are often related to continuous time systems \cite{Sun}. Sometimes, suitable symmetry conditions have been used, in order to simplify the study of such systems. In particular $Z_2$ symmetry, that is orbital symmetry with respect to one axis, has proved to be useful in approaching similar problems \cite{YH}.   \\
\indent Most of the results obtained for continuous time dynamical systems in the plane are concerned with the classical Li\'enard system (\ref{syslie})  and its generalizations, such as
\begin{equation}\label{sysGG}
\dot x = \xi(x)\bigg[ \varphi(y) - F(x)  \bigg], \qquad \dot y = -\zeta(y)g(x). 
\end{equation}
Such a class of systems also contain Lotka-Volterra systems and systems equivalent to Rayleigh equation (\ref{sysray}) as special cases. 

Even if the systems (\ref{sysGG}) reach a high level of generality, compared to Van der Pol system, 
 \\
\indent Even if the systems (\ref{sysGG}) reach a high level of generality, compared to Van der Pol system, 
$$
\dot x = y - \epsilon \left(\frac{x^3}3-x \right)   , \qquad \dot y = - x . 
$$
the first one to be  investigated  in relation to existence and uniqueness of limit cycles, an evident limitation is given by the fact that the variables $x$ and $y$ appear separately, so that mixed terms are products of single-variable functions. Since models displaying a different combination of variables do exist, different methods are desirable, in particular in absence of symmetry conditions.  \\
\indent  A recent result \cite{CRV} is concerned with systems equivalent to
\begin{equation}\label{equaCRV} 
\ddot x +  \sum_{k=0}^{N}f_{2k+1}(x){\dot x}^{2k+1} + x = 0,
\end{equation} 
with $f_{2k+1}(x)$ increasing  for $x > 0$, decreasing for $x < 0$, $k=0, \dots, N$.
On the other hand, there exist classes of second order models which are not covered by previous results. This is the case of a model developed in \cite{ETA} to describe the pedestrian's walk, which leads to the equation
\begin{equation}\label{ETA}
\ddot x + \epsilon \dot x (x^2 + x \dot x + {\dot x}^2 -1) + x = 0, \qquad \epsilon >0.
\end{equation}
Such an equation can be considered as a special case of a more general class of equations,
\begin{equation}\label{equaphi}
\ddot x + \dot x \phi(x,\dot x) + g(x) = 0.
\end{equation}
In this paper we study the class (\ref{equaphi}), assuming $ \phi(x,\dot x)$ to have strictly star-shaped level sets and $xg(x) > 0$ for $x \neq 0$. We prove a uniqueness result for limit cycles, and, under suitable additional assumptions, we show that a limit cycle exists and attracts every non-constant solution. Since $ \epsilon \dot x (x^2 + x \dot x + {\dot x}^2 -1)$ has strictly star-shaped level sets, the model introduced in \cite{ETA} has a unique limit cycle, attracting every non-constant solution. \\
\indent The result we present here is as well applicable to several equations of Li\'eanrd and Rayleigh type, in particular when they have a non-linear $g(x)$. \\
\indent This paper is organized as follows. In section 1 we study the equation (\ref{equaphi}), assuming $g(x)$ to be linear. We first prove the uniqueness theorem. The main tools applied here  is a uniqueness result proved in \cite{GS}. Then we introduce some mild additional hypotheses on the sign of $\phi(x,\dot x)$, under which the unique limit cycle attracts every non-constant solution. Then, in section 2, we assume $g(x)$ to be non-linear. We reduce the study of such a case to that of the linear $g(x)$, by  means of Conti-Filippov transformation \cite{SC}. The structure of section 2 is very similar to that of section 1, the main difference being the derivation of the condition on $ \phi(x,\dot x)$ which implies the strict star-shapedness property for the transformed system. 

\section{Linear $g(x)$}

Let $\Omega \subset \R^2$ be a star-shaped set. We denote partial derivatives by subscripts, i. e. $\phi_x$ is the derivative of $\phi$ w. r. to $x$, etc..
We say that a function $\phi \in C^1( \Omega, \R) $ is  {\it star-shaped} if $(x,y) \cdot \nabla \phi = x \phi_x + y \phi_y$ does not change sign. We say that $\phi$ is 
{\it strictly star-shaped} if  $(x,y) \cdot \nabla \phi \neq 0$, except at the origin $O=(0,0)$.  We say that $\gamma(t)$ is {\it positively bounded} if the semi-orbit $\gamma^+ = \{\gamma(t), \quad t \geq 0\}$ is contained in a bounded set. Similarly for the negative boundedness. We say that an orbit is an {\it open unbounded orbit} if it is both positively and negatively unbounded.    We say that a set $X$ is {\it invariant} if every orbit starting at a point of $X$ is entirely contained in $X$. For other definitions related to dynamical systems, we refere to \cite{BS} We call {\it ray} a half-line having origin at the point $(0,0)$.  

In this section we are concerned with the equation 
\begin{equation}\label{equaphix}
\ddot x + \dot x \phi(x,\dot x) + k x = 0, \qquad k \in \R, \qquad k > 0.
\end{equation}
Without loss of generality, possibly performing a time rescaling, we may restrict to the case $k=1$.
Let us consider a system equivalent to the equation (\ref{equaphix}), for $k=1$,
\begin{equation}\label{sysphi}
\dot x = y \qquad \dot y = -x - y \phi(x,y).
\end{equation}
We denote by $\gamma(t,x^*,y^*)$ the unique  solution of the system (\ref{sysphi}) such that $\gamma(0,x^*,y^*) = (x^*,y^*)$. We first consider a sufficient condition for limit cycles' uniqueness. We set
$$
A(x,y) = y \dot x - x \dot y = y^2 + x^ 2 +xy \phi(x,y).
$$
The sign of $A(x,y)$ is opposite to that of the angular speed of the solutions of (\ref{sysphi}).  Our uniqueness result comes from theorem 2 of \cite{GS}, in the form of corollary 6.

\begin{theorem}\label{teorema} Let $\phi\in( \R^2, \R^2)$ be a \sss function. Then (\ref{sysphi}) has at most one limit cycle.
\end{theorem}
{\it Proof.}
Without loss of generality, one may assume that, for $(x,y) \neq (0,0)$,
$x \phi_x + y  \phi_y > 0$.
The proof can be performed analogously for the opposite inequality.  \\
\indent We claim that $\nabla A(x,y)$ does not vanish on the set $A_0 = \{ (x,y) : A(x,y) =0  \} \setminus \{  (0,0) \}$. In fact,
$\nabla A$ and $A$ vanish simultaneously at $(x,y)$ if and only if
\begin{equation} \label{nabla}
\left\{  
\begin{array}{cc}  
2x + y\phi +xy \phi_x & = 0 \\ 
2y + x\phi + xy \phi_y & = 0 \\
x^ 2 + y^2 +xy \phi & = 0
\end{array} \right. .
\end{equation}
Multiplying the first equation by $y$, the second one by $x$ and re-ordering terms yields
$$
\left\{  
\begin{array}{cc} 
xy\phi & = - x^2 (2 + y\phi_x) \\ 
xy\phi & = - y^2 (2 + x\phi_y) \\ 
xy \phi & = - x^ 2 - y^2 
\end{array} \right. .
$$
Multiplying the third equation by 2 and summing with  the first two equations yields
$$
xy ( x \phi_x + y \phi_y ) = 0.
$$
Since, by hypothesis,  $x \phi_x + y \phi_y \neq 0$ except at $O$, one has $xy = 0$. If $x=0$, then by the third equation in (\ref{nabla}) one has $y=0$. 
Similarly, if  $y=0$.  \\
\indent This shows that, at every point, $A_0$ is locally a graph. Additionally, every ray $\{ (t\cos \theta, t\sin \theta), t > 0\}$, meets $A_0$ at most at a point. In fact, for $xy \neq 0$, one has, 
$$
A(t\cos \theta, t\sin \theta) = 0 \iff   \phi(t\cos \theta, t\sin \theta) = \frac{1}{ \cos \theta\sin \theta}. 
$$  
The condition $x \phi_x + y  \phi_y > 0$ implies that $\phi$ is an increasing function of $t$ on every ray. Hence on every ray not contained in an axis there exists at most one $t$ such that $\phi(t\cos \theta, t\sin \theta) = \frac{1}{ \cos \theta\sin \theta}$.  As for $xy = 0$, $A$ vanishes only at $O$.

Moreover, working as above, one can show that $A_0$ has at a single point in common with the axes, $O$. \\
\indent  The radial derivative $A_r$ of $A$ is given by
\begin{equation} \label{radder}
A_r = \frac {x A_x + y A_y}r  = \frac 1r\bigg(2 A + xy (x \phi_x + y\phi_y) \bigg) .
\end{equation}
Let $(x^*,y^*)$ be a point of the first orthant, i. e.  $x^*>0, y^* > 0$. If $A(x^*,y^*) \geq 0$, then $A_r >0$ at $(x^*,y^*)$ and at every point $(rx^*,ry^*)$ with $r > 1$, hence $A$ is strictly increasing on the half-line $\{ (rx^*,ry^*) : r > 1\}$. Now, let $(x^*,y^*)$ be a point of the second orthant, i. e.   $x^*>0, y^* < 0$. If $A(x^*,y^*) <0$, then $A_r <0$ at $(x^*,y^*)$ and at every point $(rx^*,ry^*)$ with $r > 1$, hence $A$ is strictly decreasing on the half-line $\{ (rx^*,ry^*) : r > 1\}$. The same argument allows to prove that in the third orthant $A$ behaves as in the first one, and in the fourth orthant $A$ behaves as in the second one.  \\
\indent Assume, by absurd, two distinct limit cycles to exist. The system (\ref{sysphi}) has a unique critical point, hence they are concentric. Let $\mu_1$ be the inner one, $\mu_2$ be the external one. Let $D$ be the annular region bounded by $\mu_1$ and $\mu_2$. We claim that $A(x,y) > 0$ in $D$. We prove it by proving that, for every orbit $\gamma$ contained in $D$, $A(\gamma(t)) > 0$. Let us observe that every orbit in $D$ has to meet every semi-axis, otherwise its positive limit set would  contain a critical point different from $O$. On every semi-axis one has $A(x,y) > 0$.
Assume first, by absurd, $A(\gamma(t))$ to change sign. Then there exist $t_1 < t_2$ such that $A(\gamma(t_1)) > 0$, $A(\gamma(t_2)) < 0$, and and $\gamma(t_i)$, $i = 1,2$ are on the same ray. Assume $\gamma(t_i)$, $i = 1,2$ to be in the first orthant. Two cases can occur: either  $|\gamma(t_1)| < |\gamma(t_2)|$ or  $|\gamma(t_1)| > |\gamma(t_2)|$. The former, $|\gamma(t_1)| < |\gamma(t_2)|$, contradicts the fact that $A$ is radially increasing in the first orthant, hence one has $|\gamma(t_1)| > |\gamma(t_2)|$. The orbit $\gamma(t_1)$ crosses the segment $\Sigma = \{r \gamma(t_1), 0 < r <1\}$, going towards the positive $y$-semi-axis. Let $G$ be the region bounded by the positive $y$-semi-axis, the ray $\{r \gamma(t_1), r > 0\}$ and the portions of $\mu_1$, $\mu_2$ meeting the $y$-axis and such a ray.  The orbit $\gamma$ cannot remain in $G$, since in that case $G$ would contain a critical point different from $O$. Also, $\gamma$ cannot leave $G$ crossing the positive $y$-semi-axis, because $A(x,y) > 0$ on such an axis. Hence $\gamma$ leaves $G$ passing again through the segment $\Sigma$. That implies the existence of $t_3 > t_2$, such that $\gamma(t_3)$ lies on the ray  $\{ r \gamma(t_1), 0 < r \}$. Again, one cannot have $|\gamma(t_3)| < |\gamma(t_2)|$, since $A(\gamma(t_3)) > 0$ implies $A$ increasing on the half-line $r \gamma(t_3)), r > 1$, hence one has $|\gamma(t_3)| > |\gamma(t_2)|$. Also, one cannot have $|\gamma(t_3)| < |\gamma(t_1)|$, otherwise $\gamma$ would enter a positively invariant region, bounded by the curve $\gamma(t)$, for $t_1 \leq t \leq t_3$,  and by  the segment with extrema $\gamma(t_1)$, $\gamma(t_3)$, hence there would exist a critical point different from $O$. As a consequence, one has $|\gamma(t_3)| > |\gamma(t_1)|$.
Since $A(x,y) >0$ on the segment joining $\gamma(t_1)$ and $\gamma(t_3)$, such a segment, with the portion of orbit joining $\gamma(t_1)$ and $\gamma(t_3)$ bounds a region which is negatively invariant for (\ref{sysphi}), hence contains a critical point different from $O$, contradiction.  \\
  This argument may be adapted to treat also the case of a ray in the second orthant, replacing the positive $y$-semi-axis with the positive $x$-semi-axis, and reversing the relative positions of the points $\gamma(t_1)$, $\gamma(t_2)$, $\gamma(t_3)$. In the other orthants one repeats the arguments of the first and second orthants, respectively.  \\
\indent Finally, assume that $A(\gamma(t)) = 0$ at some point $(x^*,y^*)$. If $(x^*,y^*)$ is interior to $D$, then on the ray  $\{ (rx^*,ry^*), r >0 \}$, there exist points interior to $D$ with $A(rx^*,ry^*) < 0$. Then we can apply  the above argument to the orbits starting at such points. If $(x^*,y^*)$ belongs to the boundary of $D$, then it is on $\mu_1$ or on $\mu_2$. Assume $(x^*,y^*) = \mu_1(t^*)$, for some $t^*$ (the argument works similarly on $\mu_2$). Since $A(\mu_1(t^*)) = 0$, $\mu_1$ is tangent to the ray $\{ (rx^*,ry^*), r > 0 \}$. On the other hand, $\displaystyle A_r =  \frac {xy (x \phi_x + y\phi_y) } r
\neq 0$ at $(x^*,y^*) $, hence $A_0$ is neither tangent to the ray $\{ (rx^*,ry^*), r > 0\}$, nor to $\mu_1$ at $\mu_1(t^*)$. This implies that $A_0$ and $\mu_1$ are transversal at $(x^*,y^*)$, so that a portion of $A_0$ enters $D$. Since $A_0$ separates points where $A >0$ from points where $A < 0$, also in this case there exist points interior to $D$ with $A < 0$.  \\
\indent  Now we can restrict to the annular region $D$ and divide the vector field of (\ref{sysphi}) by $A(x,y)$, as in  corollary 6 in \cite{GS}. In order to appy such a corollary, one has to  compute the expression
$$
\nu = P\left( xQ_x + yQ_y \right) - Q \left( xP_x + yP_y\right) ,
$$
where $P$ and $Q$ are the components of the considered vector field. For system (\ref{sysphi}), one has
$$
\nu =
y \left(-x - xy \phi_x - y\phi - y ^2 \phi_y \right) - \left(  -x - y \phi \right) y  =
$$
$$
-y^2 \left( x \phi_x + y  \phi_y \right) \leq 0.
$$
The function $\nu$ vanishes only for $y=0$.   For both cycles one has:
$$
\int_{0}^{T_i} \nu(\mu_i(t) )dt < 0, \qquad i=1,2,
$$
where $T_i$ is the period of $\mu_i$, $i=1,2$. Hence both cycles, by theorem 1 in \cite{GS}, are attractive. Let $A_1$ be the region of attraction of $\mu_1$.  $A_1$ is bounded, because it is enclosed by $\mu_2$, which is not attracted to $\mu_1$. The external component of   $A_1$'s boundary is itself a cycle $\mu_3$, because (\ref{sysphi}) has just one critical point at $O$. Again,
$$
\int_0^{T_3} \nu(\mu_3(t) )dt < 0,
$$
hence $\mu_3$ is attractive, too. This contradicts the fact that the solutions of (\ref{sysphi}) starting from its inner side are attracted to $\mu_1$. Hence the system  (\ref{sysphi}) can have at most a single limit cycle.
\hfill $\clubsuit$

The angular velocity of the solutions need not be negative at every point of the plane. In fact, even a simple sysytem as that one studied in (\ref{ETA}) has negative angular velocity only in a proper subset of the plane. In Figure 1 we have plotted some orbits of the equation \ref{ETA}  tending at the limit cycle, together with the two components of the curve $A_0$. The orbits cross $A_0$ at the points where their angular velocity changes sign. 

\begin{figure}[h!]
  \caption{$\ddot x + \dot x (x^2 + x \dot x + {\dot x}^2 -1) + x = 0$   }
  \centering
    \includegraphics[width=0.5\textwidth]{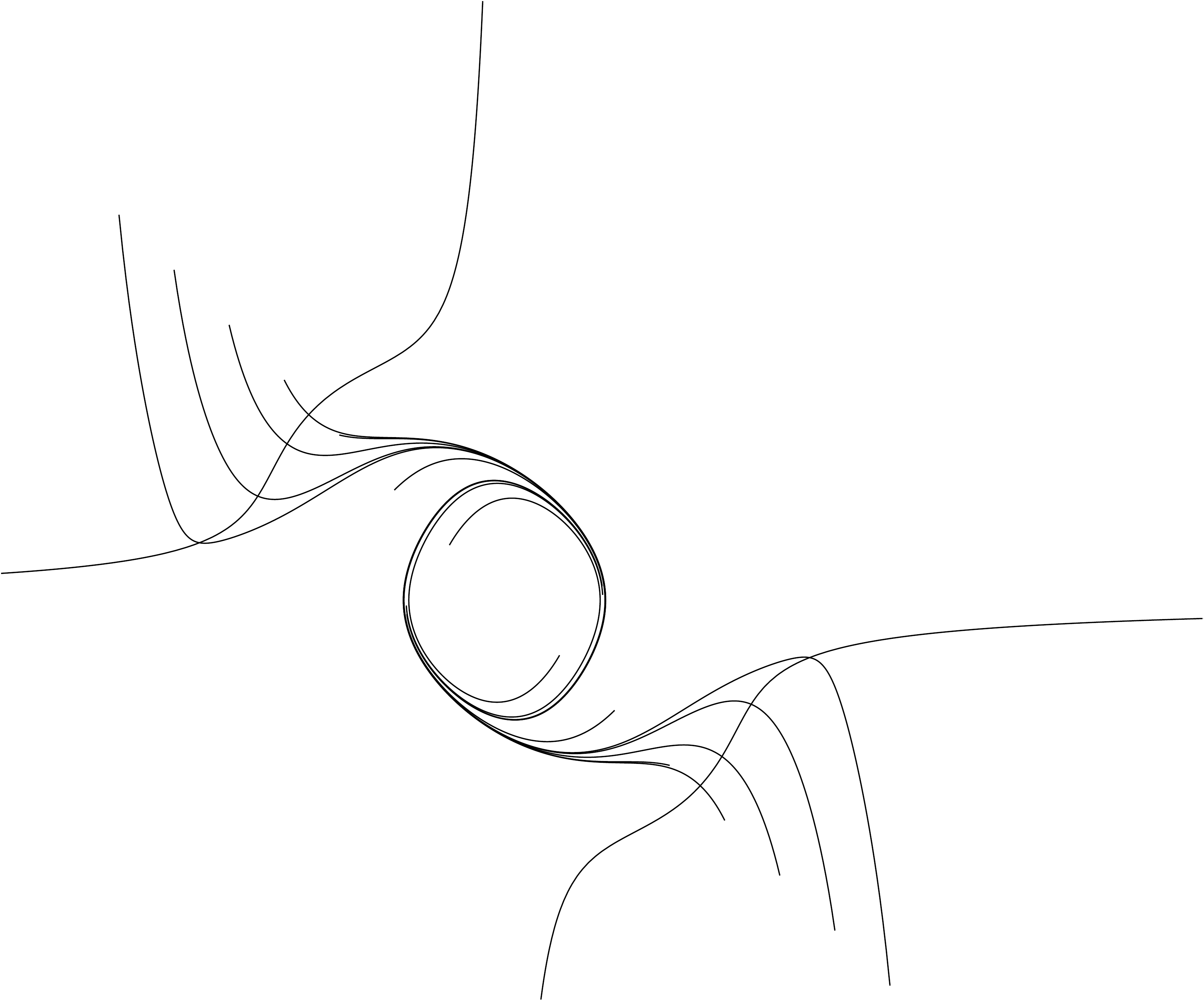}
\end{figure}

In the example of Figure 1 we have chosen $\epsilon  = 1$. In general, the system studied in \cite{ETA} has just one limit cycle, for $\epsilon  > 0$. In fact, in this case one has 
$$
\nu = x \phi_x + y  \phi_y = 2 \epsilon y^2 (x^2 + xy + y^2) > 0  \quad {\rm for } \quad  y \neq 0.
$$
It should be noted that even if the proof is essentially based on a stability argument, the divergence cannot be used in order to replace the function $\nu$. In fact, the divergence of system (\ref{sysphi}) is
$$
{\rm div } \bigg(y , -x - y \phi(x,y) \bigg)= - \phi -  y  \phi_y, 
$$
which does not have constant sign, under our assumptions. Moreover, the divergence cannot have constant sign in presence of a repelling critical point and an attracting cycle.

\indent Now we care about the existence of limit cycles.  Let us denote by $D_r$ the disk $ \{ (x,y) : x^2 + y^2 \leq r^2 \} $, and by $\partial D_r$ its boundary $\{ (x,y) : x^2 + y^2 = r^2 \} $.  Let us consider the function $V(x,y) =   \frac 12 \left( x^2 + y^2 \right)$. Its derivative along the solutions of (\ref{sysphi}) is 
$$
\dot V(x,y) = -y^2 \phi(x,y).
$$

\begin{lemma}\label{lemma} Let $U$ be a bounded set, with $\sigma := \sup \{  \sqrt{x^2 + y^2}, (x,y) \in U \}$. 
If $\phi(x,y) \geq 0$ out of $U$, and $\phi(x,y)$ does not vanish identically on any $\partial D_r$, for $r > \sigma$, then every $\gamma(t)$ definitely enters the disk $D_\sigma$ and does not leave it.
\end{lemma}
{\it Proof.}
The level curves of $V(x,y)$ are circumferences. For every $r \geq \sigma$, the disk $D_r $ contains $U$.  Since $\dot V(x,y) = - y^2 \phi(x,y) \leq 0$ on its boundary, such a disk is positively invariant. 
Let $\gamma$ be an orbit with a point $\gamma(t^*)$ such that $d^* = {\rm dist} (\gamma(t^*),O) > \sigma$. Then $\gamma$ does not leave the disk $D_{d^*}$, hence it is positively bounded.
Moreover $\gamma(t)$ cannot be definitely contained in $\partial D_r$, for any $r > \sigma$, since $\dot V(x,y)$ does not vanish identically on any $\partial D_r$, for $r > \sigma$. Now, assume by absurd that $\gamma(t)$ does not intersect  $B_\sigma$. Then its positive limit set is a cycle $\overline{\gamma}(t)$, having no points in $D_\sigma$. The cycle $\overline{\gamma}(t)$ cannot cross outwards any $\partial D_r$, hence it has to be contained in $\partial D_r$, for some $r > \sigma$, contradicting the fact that $\dot V(x,y)$ does not vanish identically on any $\partial D_r$, for $r > \sigma$. Hence there exists  $t^+ > t^*$ such that $\gamma(t^+) \in D_\sigma$. Then, for every  $t > t^+$, one has $\gamma(t) \in D_\sigma$, because $\dot V(x,y) \leq 0$ on $B_\sigma$.
\hfill$\clubsuit$

Collecting the results of the above statements, we may state a theorem of existence and uniqueness for limit cycles of a class of second order equations. We say that an equilibrium point $O$ is {\it negatively asymptotically stable} if it is asymptotically stable for the system obtained by reversing the time direction.

\begin{theorem}\label{teorema2} If the hypotheses of theorem \ref{teorema} and lemma \ref{lemma} hold, and $\phi(0,0) < 0$, then the system (\ref{sysphi}) has exactly one limit cycle, which attracts every non-constant solution.
\end{theorem}
{\it Proof.} By the above lemma, all the solutions are definitely contained in $D_\sigma$. The condition $\phi(0,0) < 0$  implies by continuity $\phi(x,y) < 0$ in a neighbourhood $N_O$ of $O$. This gives  the negative asymptotic stability of $O$ by Lasalle's invariance principle \cite{V}, since $\dot V(x,y) \geq 0$ in $N_O$, and the set $\{\dot V(x,y) = 0\} \cap N_O = \{y = 0\}   \cap N_O$ does not contain any positive semi-orbit. The system has just one critical point at $O$, hence by Poincar\'e-Bendixson theorem there exist a limit cycle. By  theorem \ref{teorema}, such a limit cycle is unique.
\hfill$\clubsuit$

This proves that every non-constant solution to the equation (\ref{ETA}) studied in  \cite{ETA} is attracted to the unique limit cycle. 

We can produce more complex systems with such a property. Let us set
$$
\phi(x,y) = -M +\sum_{k=1}^n H_{2k}(x,y),
$$
with $ H_{2k}(x,y)$ is a homogeneous function of degree $2k$, positive except at $O$, $M$ is a positive constant. Then, by Euler's identity,  one has
$$
\nu = \sum_{k=1}^n \left( x\hkx + y\hky \right) = \sum_{k=0}^n 2kH_{2k}(x,y) > 0  \quad {\rm for } \quad  (x,y) \neq (0,0).
$$
If $\phi(x,y)$ does not vanish identically on any $\partial D_r$, for instance if $H_{2k}(x,y) = (x^2 +xy + y^2)^k$, then the corresponding system (\ref{sysphi}) has a unique limit cycle.
In general, it is not necessary to assume the positiveness of all of the homogeneous functions $H_{2k}(x,y)$, as the following example shows. Let us set $Q(x,y) = x^2+xy+y^2$. Then take
$$
\phi(x,y) =  -1 + Q  - Q^2 + Q^3 .
$$
One has
$$
\nu = x \phi_x + y  \phi_y = 2Q - 4Q^2 +6Q^3 = Q(2  - 4Q  +6Q^2)  .
$$
The discriminant of the quadratic polynomial $2  - 4Q  +6Q^2$ is $ -32 < 0$ hence $\nu >0$ everywhere but at $O$. Moreover, $\phi(x,y)$ does not vanish identically on any circumference, hence  the corresponding system (\ref{sysphi}) has a unique limit cycle.

\section{Non-linear $g(x)$ }

Even if the equation we consider in this section are of a more general type, we actually derive our result from that of the previous section, so that we can consider what follows a corollary of the previous result. Let us consider the equation
\begin{equation}\label{equaphig}
\ddot x + \dot x \Phi(x,\dot x) + g(x) = 0 .
\end{equation}
We assume that  $ xg(x) > 0 {\rm\ for \ } x \neq 0$, $g \in C^1(\R,\R)$, $g'(0) \neq 0$,
We could consider equations defined in smaller subset of the plane, without essential changes. The main tools is the so-called Conti-Filippov transformation, which acts on the equivalent system
\begin{equation}\label{sysphig}
\dot x = y, \qquad  \dot y = - g(x) - y \Phi(x,y),
\end{equation}
in such a way to take the conservative part of the vector field into a linear one. Let us set $G(x) = \int_0^x g(s) ds$, and denote by $\sigma(x)$ the sign function, whose value is $-1$ for $x < 0$, $0$ at $0$, $1$ for $x > 0$.  Let us define the function $\alpha : \R \rightarrow \R$ as follows:
$$
 \alpha(x) =  \sigma(x) \sqrt{2G(x)}.
$$
Then Conti-Filippov transformation is the following one,
\begin{equation}\label{transf}
(u,v) = \Lambda(x,y) = (\alpha(x),y).
 \end{equation}
Since we assume that $g \in C^1(\R,\R)$, one has $\alpha  \in C^1(\R,\R)$ . The function $u = \alpha (x)$ is invertible, due to the condition $xg(x) > 0$. Let us call $x =\beta(u)$ its inverse. The condition $g'(0) > 0$ guarantees the differentiability of $\beta(u)$ at $O$. For $x \neq 0$, that is for $u \neq 0$, one has,
\begin{equation}\label{derivate}
\alpha'(x) =  \frac{\sigma(x) g(x)}{\sqrt{2G(x)}}, \qquad \beta'(u) =  \frac{1}{\alpha'(\beta(u))} =  \frac{\sigma(x)\sqrt{2G(x)}}{g(x)} = 	\frac{u}{g(\beta(u))}.
 \end{equation}
 For $x=u=0$ one has, 
 $$
 \alpha'(0) = \sqrt{g'(0)}, \qquad  \beta'(0) =  \sqrt{\frac{1}{g'(0)}} .
 $$
Finally,
$$
\lim_{u \rightarrow 0} \frac{g(\beta(u))}{u} =  \sqrt{g'(0)} > 0.
$$
\begin{theorem}\label{teorema-n} Assume $g \in C^1(\R,\R)$, with $g'(x) > 0$ and $xg(x) > 0$ for $x\neq 0$. Let $\Phi \in C^1(\R^2, \R^2) $ satisfy
$$
 \frac {\sigma(x) \sqrt{2G(x)} }{g(x)} \left[ 2G(x) \frac{ \Phi_x(x,y)g(x) - \Phi(x,y)g'(x) }{g(x)^2} +  \Phi(x,y)  \right]
+ y \Phi_y(x,y) \neq 0.
$$
Then (\ref{sysphig}) has at most one limit cycle.
\end{theorem}
{\it Proof.}
For $u \neq 0$, the transformed system has the form
\begin{equation}\label{sysu}
\dot u = v \frac{g(\beta(u))}{u}, \qquad \dot v = - g(\beta(u)) - v \Phi(\beta(u),v).
\end{equation}
For $u = 0$, the above form is extended by continuity. In the following we consider only the case $u\neq 0 $, that is $x \neq 0$, since the case $u = 0 = x$ is obtained by continuity. We may multiply the system (\ref{sysu}) by $ \frac{u}{g(\beta(u))}$, obtaining a new system having the same orbits as (\ref{sysu}), 
 \begin{equation}\label{sysurep}
\dot u = v, \qquad \dot v = - u - v  \frac{u\Phi(\beta(u),v)}{g(\beta(u))}.
\end{equation}
Such a system is of the type (\ref{sysphi}), so that we may apply theorem \ref{teorema} to get uniqueness of solutions. This reduces to require the strict star-shapedness of the function
$$
 \frac{u\Phi(\beta(u),v)}{g(\beta(u))},
$$
that is,
$$
u\left[  \frac{u\Phi(\beta(u),v)}{g(\beta(u))}  \right]_u+ v \Phi_v(\beta(u),v) > 0.
$$
Since $\beta(u) = x$, $v = y$, the second term in the above sum is just $y \Phi_y(x,y)$. As for the the first one, one has
$$
\left[  \frac{u\Phi(\beta(u),v)}{g(\beta(u))} \right]_u = $$ $$
 \frac{\bigg( u\Phi_u(\beta(u),v)\beta'(u) +  \Phi(\beta(u),v) \bigg)g(\beta(u)) - u\Phi(\beta(u),v)g'(\beta(u)) \beta'(u)}{g(\beta(u))^2}.
$$
Replacing $\beta(u)$ with $x$ and applying the formulae (\ref{derivate}) one has 
$$
\left[  \frac{u\Phi(\beta(u),v)}{g(\beta(u))} \right]_u = $$ $$
 \frac{\sigma(x)\sqrt{2G(x)}}{g(x)}.
 \frac{\bigg( \sigma(x) \sqrt{2G(x)}\Phi_x(x,y) \bigg)g(x) - \sigma(x) \sqrt{2G(x)}\Phi(x,y)g'(x) }{g(x)^2} + 
 \frac{\Phi(x,y) }{g(x)}.
$$
Since $\sigma(x)^2 = 1$ eveywhere but at $0$, the above formula reduces to
$$
\left[  \frac{u\Phi(\beta(u),v)}{g(\beta(u))} \right]_u =
 \frac{2G(x)}{g(x)}   \frac{ \Phi_x(x,y)g(x) - \Phi(x,y)g'(x) }{g(x)^2} +   \frac{\Phi(x,y) }{g(x)} .
$$
Concluding, one has
$$
u\left[  \frac{u\Phi(\beta(u),v)}{g(\beta(u))}  \right]_u = 
 \frac {\sigma(x) \sqrt{2G(x)} }{g(x)} \left[ 2G(x) \frac{ \Phi_x(x,y)g(x) - \Phi(x,y)g'(x) }{g(x)^2} +  \Phi(x,y)  \right].
$$
Hence, the star-shapedness conditions reduces to
$$
 \frac {\sigma(x) \sqrt{2G(x)} }{g(x)} \left[ 2G(x) \frac{ \Phi_x(x,y)g(x) - \Phi(x,y)g'(x) }{g(x)^2} +  \Phi(x,y)  \right]
+ y \Phi_y(x,y) \neq 0.
$$
\hfill$\clubsuit$

If $g(x) = x$, then $\displaystyle{ G(x) =  \frac{x^2}{2} }$ and $\displaystyle{  \frac {\sigma(x) \sqrt{2G(x)} }{g(x)} = 1 }$, for $x \neq 0$. In this case the star-shapedness condition just reduces to what considered in the previous section, since
$$
 2G(x)  \frac{ \Phi_x(x,y)g(x) - \Phi(x,y)g'(x) }{g(x)^2} +  \Phi(x,y)  =   x \Phi_x(x,y) .
$$

Now we prove the non-linear analogous of lemma \ref{lemma}. Let us set
$$
E(x,y) = G(x) + \frac {y^2}{2}
$$
For $r > 0$, we set $\Delta_r = \{ (x,y) : 2 E(x,y) < r^2 \}$ and $\partial \Delta_r = \{ (x,y) : 2 E(x,y)  = r^2 \}$

\begin{lemma}\label{lemma-n} Let $U$ be a bounded set, with $\sigma := \sup \{   \sqrt{ 2 E(x,y) }, (x,y) \in U \}$. 
If $\Phi(x,y) \geq 0$ out of $U$, and $\Phi(x,y)$ does not vanish identically on any $\partial \Delta_r$, for $r > \sigma$, then every $\gamma(t)$ definitely enters the set $\Delta_\sigma$ and does not leave it.
\end{lemma}
{\it Proof.} Performing Conti-Filippov transformation, the sets $ \Delta_r$ are taken into the sets $D_r$, as well as the boundaries $\partial \Delta_r$ are taken into the boundaries $\partial D_r$. The function $\Phi(x,y)$ does not vanish identically on any $\partial \Delta_r$ if and only if the function $\phi(u,v) = \Phi(\beta(u),v)$ does not vanish identically on any $\partial D_r$.
Then one can apply lemma  \ref{lemma} to the system (\ref{sysurep}). In fact, the derivative of the Liapunov function $V(u,v) = \frac 12 (u^2+v^2)$ along the solutions of  (\ref{sysurep}) is just 
$$ 
\dot V (u,v) =  - v^2  \frac{u\Phi(\beta(u),v)} {g(\beta(u))} .
$$
The function $ \frac{u} {g(\beta(u))}$ is positive for $u \neq 0$, hence the hypotheses of lemma \ref{lemma} are satisfied by the system (\ref{sysurep}). As a consequence, the conclusions of lemma  \ref{lemma} hold for the system (\ref{sysurep}), and applyng the inverse transformation $\Lambda^{-1}$ one obtains the thesis.
\hfill$\clubsuit$

Now we can conclude proving the analogue of theorem \ref{teorema2} for the equation with a non-linear  $g(x)$.

\begin{theorem}\label{teorema2-n}  If the hypotheses of theorem \ref{teorema-n} and lemma \ref{lemma-n} hold, and $\phi(0,0) < 0$, then the system (\ref{sysphig}) has exactly one limit cycle, which attracts every non-constant solution.
\end{theorem}
{\it Proof.} As the proof of theorem \ref{teorema2}, replacing the Liapunov function $V(x,y)$ with the Liapunov function $E(x,y)$.
\hfill$\clubsuit$

\bigskip

{\bf  \large Acknowledgements}

The author would like to thank dr. S. Erlicher for raising the problem, and proff. T. Carletti and G. Villari for reading a previous version of this paper. 

This paper has been partially supported by the GNAMPA 2009 project \lq\lq Studio delle traiettorie di equazioni differenziali ordinarie\rq\rq.

\end{document}